\documentclass[12pt]{amsart}
\usepackage{fullpage}
\usepackage{amssymb}
\usepackage{graphicx}

\newtheorem{theorem}{Theorem!!!}[section]
\newtheorem{thm}[theorem]{Theorem}
\newtheorem{cor}[theorem]{Corollary}
\newtheorem{prop}[theorem]{Proposition}
\newtheorem{lem}[theorem]{Lemma}

\theoremstyle{definition}
\newtheorem{defn}{Definition}[section]

\newcommand{\thmref}[1]{Theorem \ref{#1}}

\newcommand{\comm}[1]{}

\newcommand{\BBB}[1]{{\mathbb #1}}
\newcommand{\ra}{\rightarrow}

\newcommand{\dist}{\operatorname{dist}}

\newcommand{\eps}{\epsilon}

\newcommand{\Om}{\Omega}
\newcommand{\ze}{\zeta}

\newcommand{\tl}{\tilde}
\newcommand{\NN}{{\BBB N}}
\newcommand{\ZZ}{{\BBB Z}}

\newcommand{\JJ}{{\BBB J}}
\newcommand{\RR}{{\BBB R}}

\newcommand{\CC}{{\BBB C}}
\newcommand{\DD}{{\BBB D}}

\newcommand{\cC}{{\mathcal C}}
\newcommand{\cF}{{\mathcal F}}

\begin{document}

\title[On computability of Julia sets]{On computability of Julia sets: answers to questions
of Milnor and Shub}

\begin{abstract}
In this note we give answers to questions posed to us by J.~Milnor and M.~Shub,
which shed further light on the structure of non-computable Julia sets.
\end{abstract}

\author{M. Braverman}
\author{M. Yampolsky}
\date{\today}
\maketitle

\section{Introduction}
\subsection*{Computability of real sets.}
The reader is directed to \cite{BY} for a more detailed
 discussion of the notion of computability of
subsets of $\RR^n$ as applied, in particular, to Julia sets. We recall the 
principal definitions here. The exposition below uses the concept of a {\it Turing
Machine}. This is a standard model for a computer program employed by computer
scientists. Readers unfamiliar with this concept should think instead of an
algorithm written in their favorite programming language. These concepts are
known to be equivalent.

Denote by $\DD$ the set of the {\it dyadic rationals}, that is, rationals of 
the form $\frac{p}{2^m}$. 
%Rationals in $\DD$ can be easily represented as binary numbers. 
We say that $\phi: \NN \rightarrow \DD$ is an {\it oracle}
for a real number $x$, if $| x - \phi(n)|<2^{-n}$ for all $n \in \NN$. 
In other words, $\phi$ provides a good dyadic approximation for $x$. 
We say that a Turing Machine (further abbreviated as TM) $M^{\phi}$ is 
an {\it oracle machine}, if at every step of 
the computation M is allowed to query the value $\phi(n)$ for any $n$.
This definition allows us to define the computability of real functions 
on compact sets. 

\begin{defn}
\label{funcomp}
We say that a function $f:[a,b] \rightarrow [c,d]$ is computable, if 
there exists an oracle TM $M^{\phi} (m)$ such that if $\phi$ is an oracle 
for $x \in [a,b]$, then on input $m$, $M^{\phi}$ outputs a $y \in \DD$ 
such that $| y - f(x)|<2^{-m}$.
\end{defn}

%\noindent
%To understand this definition better, the reader without a Computer Science
%background should think of a computer program with an instruction
%$$\text{ READ real number }x\text{ WITH PRECISION }n(m).$$
%On the execution of this command, a dyadic rational $d$ is input from the 
%keyboard. This number must not differ from $x$ by more than  $2^{-n(m)}$ 
%(but otherwise can be arbitrary). The algorithm then outputs $f(x)$ to precision
%$2^{-n}$.

\noindent
It is worthwhile to note why the oracle mechanism is introduced. 
There are only countably many possible algorithms, and consequently only
countably many {\it computable} real numbers which such algorithms can
encode. Therefore, one wants to separate the hardness of encoding the
real number $x$ from the hardness of computing the value of the function
$f(x)$, having the access to the value of $x$.

\noindent
Let $K \subset \RR^k$ be a compact set. 
We say that a TM M computes the set $K$ if it approximates $K$ in the {\it 
Hausdorff metric}. Recall that the Hausdorff metric is a metric on 
compact subsets of $\RR^n$ defined by 
\begin{equation}
\label{hausdorff metric}
d_H ( X, Y) =  \inf \{\epsilon > 0 \;|\; X \subset U_{\epsilon} 
(Y)~~\mbox{and}~~  Y \subset U_{\epsilon}(X)\},
\end{equation}

\noindent
where $U_\eps(S)$ is defined as the union of the set of $\eps$-balls with centers in $S$.

\noindent
 We introduce
a class $\cC$ of sets which is dense in metric $d_H$ among the compact 
sets and which has a natural correspondence to binary strings. 
Namely $\cC$ is the set of finite unions of dyadic balls:
$$
\cC= \left\{ \bigcup_{i=1}^n \overline{B(d_i, r_i)}~|~~\mbox{where}~~d_i, 
r_i 
\in \DD \right\}.
$$
Members of $\cC$ can be encoded as binary strings in a natural way. 

\noindent
We now define the notion of computability of subsets of $\RR^n$ (see \cite{Wei}, and
 also \cite{WeiPaper}).

\begin{defn}
\label{setcomp}
We say that a compact set $K\subset \RR^k$ is computable, if there exists a
TM $M(d,n)$, where $d\in\DD$, $n\in \NN$ which outputs a value $1$ if $\dist(d,K)<2^{-n}$,
the value $0$ if $\dist(d,K)>2\cdot 2^{-n}$, and in the ``in-between'' case it halts and
outputs either $0$ or $1$. 

In other words, it
computes, in the classical sense,  a function from the family $\cF_K$
of functions of the form
\begin{equation}
\label{defcomp}
f(d,n)=\left\{
\begin{array}{ll}
0,& \text{if }\dist(d,K)>2\cdot 2^{-n}\\
1,& \text{if }\dist(d,K)<2^{-n}\\
0\text{ or }1,&\text{otherwise}
\end{array}
\right.
\end{equation} 
\end{defn}

\begin{thm}
\label{equiv-defn}
For a compact $K \subset \RR^k$ the following are equivalent:

(1) $K$ is computable as per definition \ref{setcomp},

(2) there exists a TM 
$M(m)$, such that on input $m$, $M(m)$ outputs an encoding of $C_m \in 
\cC$ such that $d_H (K, C_m) < 2^{-m}$ (global computability),

(3) the {\em distance function} $d_K (x) = \inf \{ |x-y|~~|~~y\in K \}$ 
is 
computable as per definition \ref{funcomp}.
\end{thm}

Note that in the case $k=2$  computability means that $K$ can be drawn on a computer screen 
with arbitrarily good precision (if we imagine the screen as a lattice of pixels).

In the present paper we are interested in   questions concerning 
the computability of the Julia set $J_c = J(f_c) = J(z^2+c)$.
Since there 
are uncountably many possible parameter values for $c$, we cannot 
expect for each $c$ to have a machine $M$ such that $M$ computes 
$J_c$ (recall that there are countably many TMs). On the other hand, 
it is reasonable to want $M$ to compute $J_c$ with an oracle 
access to $c$. Define the function $J: \CC \rightarrow K^{*}$ ($K^{*}$ is 
the set of all compact subsets of $\CC$) by $J(c)=J(f_c)$. 
In a complete analogy to Definition \ref{funcomp} we can define

\begin{defn}
\label{funcomp2}
We say that a function $\kappa:S \rightarrow K^{*}$ for some bounded set $S$ is 
computable, if there exits an oracle TM $M^{\phi}(d,n)$, 
where $\phi$ is an oracle for $x\in S$, which computes a function (\ref{defcomp}) of the family
$\cF_{\kappa(x)}$.

Equivalently, there exists an oracle TM $M^\phi(m)$ with $\phi$ again representing
$x\in S$ such that 
 on input $m$, $M^{\phi}$ outputs a $C \in \cC$ 
such that $d_H (C, \kappa(x))<2^{-m}$.
\end{defn}

\noindent
In the case of Julia sets:

\begin{defn}
\label{Jcomp}
We say that $J_c$ is computable if the function $J: d \mapsto J_d$ is 
computable on the set $\{ c \}$.
\end{defn}

\noindent
We have the following (see \cite{BBY2}):
\begin{thm}
Suppose that a TM $M^\phi$ computes the function $J$ on a set $S\subset\CC$. Then $J$ is continuous
on $S$ in Hausdorff sense.
\end{thm}

\subsection*{Previous results.}
We have demonstrated in \cite{BY}:

\begin{thm}
\label{non-comput}
There exists a parameter value $c\in \CC$ such that the Julia set of the
quadratic polynomial $f_c(z)=z^2+c$ is not computable.
\end{thm}

\noindent
The quadratic polynomials in \thmref{non-comput} possess Siegel disks.
It was further shown by I.~Binder and the authors of the present paper in \cite{BBY1}
that the absence of rotation domains, that is either Siegel disks or
Herman rings, guarantees computability of the rational Julia set.
This implies, in particular, that all Cremer quadratic Julia sets are computable --
this despite the fact that no informative high resolution images of such
sets have ever been produced. One expects, however, that such ``bad'' but
still computable examples have high algorithmic complexity, which makes
the computational cost of producing such a picture prohibitively high.

\subsection*{Two questions on computability of Julia sets.}
J. Milnor has asked us the following natural question:

\medskip
\noindent
{\it Is the filled Julia set of a quadratic polynomial always computable?}

\medskip
\noindent
In this paper we answer in the affirmative:

\begin{thm}
\label{thm:filled}
For any polynomial $p(z)$ there is an oracle Turing Machine
$M^{\phi}(n)$ that given an oracle access to the coefficients
of $p(z)$ and $n$, outputs a $2^{-n}$-approximation of 
the filled Julia set $K_{p(z)}$. 

Moreover, in the case when $p(z)=z^2+c$ is quadratic, only two 
machines suffice to compute all non-parabolic Julia sets: one 
for $c\in M$, and one for $c\notin M$. 
\end{thm}

\begin{figure}[!h]
\begin{center}
\includegraphics[angle=0, scale=0.5]{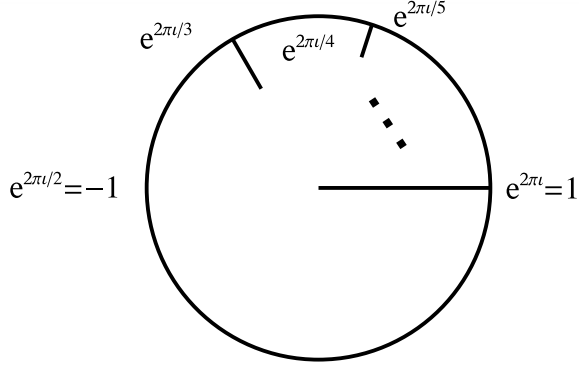}
\end{center}
\caption{Part of the picture of $\Om_t$ for $t=(0.10101\ldots)_2$}
\label{fig:Om}
\end{figure}

\noindent
This may come as a surprise, given the negative result for Julia sets. To gain some 
insight into how non-computability can be destroyed by filling in, consider the
following toy example. Let $A: \NN \ra \{0,1\}$ be any uncomputable 
predicate. Consider the set 
$$\Omega_t=
\left\{
\begin{array}{ll}
S^1\underset{k\in\NN,\;d_k=1}{\bigcup}\{r e^{2\pi i/k}|\; r\in[1-\frac{1}{k},1]\} & \mbox{for } t=(0.d_1d_2d_3\ldots)_2 \in [0,1) \\ \\
S^1\underset{k\in\NN,\;A(k)=1}{\bigcup}\{r e^{2\pi i/k}|\; r\in[1-\frac{1}{k},1]\} & \mbox{for }  t=1 
\end{array}
\right.$$
To avoid ambiguity, we always take the finite expansion for dyadic $t$'s. 
An example of a set $\Om_t$ is depicted on Figure \ref{fig:Om}.
Firstly, note that if $t\in(0,1)$ is not a computable real, then the set
$\Omega_t$ is non-computable by a TM {\it without} an oracle for $t$.
Moreover, even for a TM $M^\phi$ equipped with an oracle input for $t$, the set
$\Om_1$ is clearly non-computable.
 However, when filled, every $\Om_t$ becomes a computable set -- the unit disk.

\noindent
The question of M. Shub again has to do with fragility of
non-computability. This time, instead of filling in a 
non-computable Julia set, we will make a ``fuzzy''
picture of it by letting the parameter $c$ vary in a neighborhood.

To formalize this, consider the following definition.
Let $\JJ$ be the subset of $\CC \times \CC$ given by
$$
\JJ = \overline{\{ (z,c)~:~z\in J_c\}}. 
$$
Shub has asked us:

\medskip
\noindent
{\it Is the set $\JJ$ computable?}
\medskip
\noindent

The answer again is ``yes'':

\begin{thm}
\label{thm:main}
Let $d>0$ be any computable real. Then the set
$$\JJ\cap \CC\times \overline{B(0,d)}$$ is a computable subset of $\CC \times \CC$.
\end{thm}

\noindent
Informally, we may think of projection of $\JJ\cap \CC\times(c-\eps,c+\eps)$ 
to the first coordinate as the
picture that a computer would produce when $J_c$ itself is uncomputable.

To understand how the mechanism of non-computability is destroyed in this case,
consider again the set $\Omega_t$ for $t\in(0,1]$ as the toy model.
The set $W=\overline{\{ (z,t)~:~z\in \Omega_t,\;t\in(0,1]\}}\subset \CC\times\RR$ is  
computable even though $\Omega_t$ itself is non-computable for $t=1$. This happens
because in the closure of $W$ the ``slice" corresponding to $t=1$ is 
$$
S^1\underset{k\in\NN}{\bigcup}\{r e^{2\pi i/k}|\; r\in[1-\frac{1}{k},1]\} \supset \Om_1.
$$
This set ``masks" the computational hardness of $\Om_1$, and makes $W$ 
computable.

\section{Computability of filled Julia sets}

Our goal in this section is to prove Theorem \ref{thm:filled}. 
For a given polynomial $p(z)$ we construct a machine computing 
the corresponding filled Julia set $K_p$. We will use 
the following combinatorial information about $p$ in the construction.
Note the this information can be encoded using a finite number 
of bits. 

\begin{itemize}
\item 
Information that would allow us to compute
the non-repelling orbits of the polynomial
with an arbitrary precision. 
Note that there are at most $\deg p-1$ of them. 
Such information could, for example, consist
of the list of periods $k_i$ of such orbits;
and for each $i$ a finite collection of dyadic
balls $\{D_i^j\}_{j=1}^{k_i}$ separating the points of the 
corresponding orbit from
the other solutions of the equation $p^{k_i}(z)=z$.
This allows for an arbitrarily precise approximation
of the orbits by using an iterative root-finding algorithm
for $p^{k_i}(z)=z$ in $D_i^j$.
\item
In the case of a hyperbolic or a parabolic orbit $\ze$, a 
domain of attraction $D_{\ze}$ such that every orbit converging to 
$\ze$ eventually reaches and stays in $D_{\ze}$. In the
hyperbolic case $D_\ze$ is just a collection of discs. In 
the parabolic case, it is a collection of sectors around 
the points of $\ze$, and it can be computed with an 
arbitrarily high precision. 
\item 
In the case of a Siegel disc $D$, information that would allow us
to identify a repelling orbit $\ze_D$ in the same connected component
of $K_p$ as $D$. Such an orbit always exists, and can be identified
using a finite amount of combinatorial information.  
\end{itemize}

\subsection{Computing $K_p$}

We are given a dyadic point $d \in \DD$ and an $n\in \NN$. Our goal is 
to always terminate and output $1$ if $B(d,2^{-n}) \cap K_p \neq
\emptyset$ and to output $0$ if $B(d,2\cdot 2^{-n}) \cap K_p = \emptyset$. 
We do it by constructing five machines. They are guaranteed to 
terminate each on a different condition, always with a valid answer. 
Together they cover all the possible cases. 

\begin{lem}
\label{lem:filled}
There are five oracle machines $M_{ext}$, $M_{jul}$, $M_{hyp}$, $M_{par}$, 
$M_{sieg}$ such that 
\begin{enumerate}
\item
if $d$ is at distance $\ge \frac{4}{3}\cdot 2^{-n}$ from $K_p$, 
$M_{ext}(d,n)$ will halt and output $0$. If $d$ is 
at distance $\le 2^{-n}$ from $K_p$, $M_{ext}(d,n)$ will 
never halt;
\item
if $d$ is at distance $\le \frac{5}{3} \cdot 2^{-n}$ from $J_p$, 
$M_{jul}(d,n)$ will halt and output $1$. If $d$ is at
distance $\ge 2\cdot 2^{-n}$ from $J_p$, $M_{jul}(d,n)$  will
never halt; 
\item 
$M_{hyp}(d,n)$ halts and outputs $1$ if and only if $d$ is inside
an attracting basin of a  hyperbolic orbit of $p$;
\item
$M_{par}(d,n)$ halts and outputs $1$ if and only if $d$ is inside
an attracting basin of a parabolic orbit of $p$;
\item 
$M_{sieg}(d,n)$ halts and outputs
$1$ if the orbit of $d$ reaches a Siegel disc, and 
$d$ is at distance $\ge \frac{4}{3} \cdot 2^{-n}$ from $J_p$. 
It never halts if $d$ is at distance $\ge 2 \cdot 2^{-n}$
from $K_p$. 
\end{enumerate}
\end{lem}

\noindent
Recall the Fatou-Sullivan classification of Fatou components of
a polynomial mapping of $\hat\CC$ (see e.g. \cite{Mil}):

\begin{thm}[{\bf Fatou-Sullivan classification}]
\label{nowandering}
Every Fatou component of a polynomial mapping of $\hat\CC$ of degree at 
least two is a preimage of a periodic component. Every periodic component
is of one of the following types: the immediate basin of an attracting 
(or a super-attracting) periodic point; a component of the immediate basin
of a parabolic periodic point; a Siegel disk. 
\end{thm}

\begin{proof}[Proof of Theorem \ref{thm:filled}, given Lemma 
\ref{lem:filled}] By  Fatou-Sullivan classification 
 it is not hard to see that for each $(d,n)$ at
least one of the machines halts. Moreover, by the definition 
of the machines, they always output a valid answer whenever they
halt. Hence running the machines in parallel and returning
the output of the first machine to halt gives the algorithm 
for computing $K_p$. 
\end{proof}

It remains to prove Lemma \ref{lem:filled}.

\begin{proof} {\bf (of Lemma \ref{lem:filled})}
We give a simple construction for each of the five machines. 
\begin{enumerate}
\item
$M_{ext}$: Take a large ball $B$ such that $K_p \subset B$. 
Intuitively, we pull the ball back under $p$ to get a good 
approximation of $K_p$. Let $B_k$
be a $2^{-(n+3)}$-approximation of the set $p^{-k}(B)$. Output 
$0$ iff $B_k \cap B(d,\frac{7}{6}\cdot 2^{-n}) = \emptyset$. 
It is not hard to see that this algorithm satisfies the conditions 
on $M_{ext}$. 
\item 
$M_{jul}$: Enumerate all the repelling periodic orbits of $p$. Let $C_k$
be a $2^{-(n+3)}$-approximation of the union of the first $k$ 
orbits enumerated. Output $1$ iff $d(d,C_k) < \frac{11}{6}\cdot 2^{-n}$.
The repelling periodic orbits are all in $J_p$ and are dense in
this set. Hence the algorithm satisfies the conditions on $M_{jul}$. 
\item 
$M_{hyp}$: Let $z_k$ be a $2^{-k}$-approximation of $p^k(d)$. 
If $d$ isinside the  basin of attraction for a some orbit $\ze$, then 
$z_k$ for some $k$ will be inside $D_\ze$ for some $k$. Output $1$ 
if $z_k$ is at least $2^{-k}$-far from the boundary of $D_\ze$. 
\item 
$M_{par}$: Very similar to $M_{hyp}$. The only difference is
that now we are checking for convergence to an attracting petal 
of a parabolic orbit. 
\item 
$M_{sieg}$: This is the most interesting case. It is not 
hard to see that for each $k$, we can compute a union $E_k$ of 
dyadic balls such that 
$$
\bigcup_{i=0}^k p^i \left(B(d,\frac{4}{3} \cdot 2^{-n})\right) \subset 
E_k \subset \bigcup_{i=0}^k p^i \left(B(d, \frac{5}{3} \cdot 
2^{-n})\right).
$$
Let $c$ be the center of the Siegel disc (one of the centers, in case 
of an orbit), and let $y$ be the given periodic point in the connected 
component of $c$. We terminate and output $1$ if $E_n$ separates $c$ from 
$y$ in $\CC$ (or covers either one of them). 

If $d$ is inside the Siegel disc, then the forward images of
$B(d, \frac{4}{3} \cdot 2^{-n})$ will cover an annulus in the 
disc that will separate $c$ from the boundary of the disc, and
in particular from $y$. Hence $M_{sieg}$ will terminate and output 
$1$. 

If the distance from $d$ to $K_c$ is $\ge 2\cdot 2^{-n}$, then
$E_k \cap K_p = \emptyset$ for all $k$. In particular, $E_k$ 
cannot separate $c$ from $y$, since they are connected in $K_p$.
\end{enumerate}
\end{proof}

\subsection{The quadratic case}

In the quadratic case there is at most one non-repelling orbit. 
In the case there is a Siegel disc, the Julia set is connected, and
any repelling periodic orbit can be taken as the orbit connected to 
the Siegel disc. In fact, it is not hard to see that if we exclude 
the parabolic case, one machine suffices to take care of all the 
{\em connected} filled Julia sets. As a corollary we get:

\begin{cor}
Denote by $M$ the Mandelbrot set, and by $P$ the set of 
$c$'s for which $J_c$ is parabolic. 
The function $K: c \mapsto K_{z^2+c}$ is continuous in the 
Hausdorff metric on the set $M-P$. 
\end{cor}

\section{Computability of the set $\JJ$}

Recall that 
$$
\JJ = \overline{\{ (z,c)~:~z\in J_c\}} \subset \CC^2. 
$$
Theorem \ref{thm:main} asserts that $\JJ$ is computable. 
We prove it by showing that $\JJ$ is {\em weakly} computable. 

\begin{defn} 
\label{weakdef}
We say that a set $C$ is {\em weakly} computable if there is an oracle
Turing Machine $M^{\phi}(n)$ such that if $\phi$ represents a real number
$x$, then the output of $M^{\phi}(n)$ is
$$
M^{\phi}(n) = \left\{
\begin{array}{ll}
1 & \mbox{  if  } x \in C \\
0 & \mbox{  if  } B(x,  2^{-n}) \cap C = \emptyset \\
0 \mbox{  or  } 1 & \mbox{  otherwise  }
\end{array}
\right.
$$
\end{defn}
It has been shown that the weak definition is equivalent to the standard definition. See 
\cite{Brv05}, for example. 
We will need the following lemma. 

\begin{lem}
\label{lem:onlyhyp}
For any point $(z,c)$ in the complement of the closure $\overline{\JJ}$, 
$z$ converges to an attracting periodic orbit of $f_c: z \mapsto z^2+c$. 
\end{lem}

The proof of the lemma occupies \S\ref{sec:perturb}.

The following lemma allows us to ``cover" all points that belong to
$\JJ$.

\begin{lem}
\label{lem:alg1}
There is an algorithm $A_1 (n)$ that  on input $n$ outputs a sequence of 
dyadic points $p_1, p_2, \ldots \in \CC \times \CC$ such that 
$$
B(\JJ,2^{-(n+3)}) \subset \bigcup_{j=1}^{\infty} B(p_j, 2^{-(n+2)}) \subset B(\JJ, 2^{-(n+1)}). 
$$
\end{lem}

\begin{proof}
It is well known that repelling periodic orbits of $f_c$ are dense in $J_c$. 
Hence, the set 
$$
S_{rep} = \{(z,c)~:~z\mbox{ is in a repelling periodic orbit of }f_c\}
$$
is dense in $\JJ$. $S_{rep}$ is a union of a countable number of algebraic 
curves $S^m_{rep}$ given by the constraints 
$$
\left\{
\begin{array}{l}
f_c^m(z)=z \\
|(f_c^m)'(z)|>1
\end{array}
\right.
$$
For each $m$ we can compute a finite number of points $p_1^m, \ldots, p_{r_m}^m$ approximating 
$S^m_{rep}$ such that 
$$
B(S^m_{rep},2^{-(n+3)}) \subset \bigcup_{j=1}^{r_m} B(p_j^m, 2^{-(n+2)}) \subset B(S^m_{rep}, 2^{-(n+1)}). 
$$
We have 
$$
\overline{\JJ} = \overline{S_{rep}} = \overline{\bigcup_{m=1}^\infty S^m_{rep}}. 
$$
Hence the computable sequence $p_1^1, \ldots, p_{r_1}^1, p_1^2, \ldots, p_{r_2}^2, \ldots, p_1^m, \ldots, p_{r_m}^m,\ldots$ satisfies the conditions of the lemma. 
\end{proof}

\begin{cor}
\label{cor:alg1}
There is an oracle machine $M_1^{\phi_1, \phi_2}(n)$, where $\phi_1$ is an oracle for 
$z \in \CC$ and $\phi_2$ is an oracle for $c \in \CC$, such that $M_1^{\phi_1, \phi_2}$ always
halts whenever $d((z,c),\JJ)<2^{-(n+4)}$ and never halts if $d((z,c),\JJ)\ge 2^{-n}$.
\end{cor}

\begin{proof}
Query the oracles for a point $p \in \CC\times \CC$ such that $d(p,(z,c))<2^{-(n+4)}$. 
Then run the following loop:

\medskip
\noindent
$i \leftarrow 0$ \\
 {\bf do}\\
$~~~~$ $i \leftarrow i+1$\\
$~~~~$ generate $p_i$ using $A_1(n)$ from Lemma \ref{lem:alg1} \\
{\bf while} $d(p,p_i)>2^{-(n+2)}$   
\medskip

If $d((z,c),\JJ)<2^{-(n+4)}$, then $d(p, \JJ)<2^{-(n+3)}$, hence by Lemma 
\ref{lem:alg1} there is an $i$ such that $d(p,p_j)\le 2^{-(n+2)}$, and the 
loop terminates. If $d((z,c),\JJ)>2^{-n}$, then $d(p,\JJ)>2^{-n}-2^{-(n-4)} >1.5 \cdot 2^{-(n+1)}$. 
Hence, by Lemma \ref{lem:alg1}, $p \notin B(p_i, 2^{-(n+1)})$ for all $i$, and the loop will 
never terminate. 
\end{proof}

The following lemma allows us to exclude points outside $\overline{\JJ}$ from 
$\JJ$. 

\begin{lem}
\label{alg2}
There is an oracle machine $M_2^{\phi_1, \phi_2}$, where $\phi_1$ is an oracle for 
$z \in \CC$ and $\phi_2$ is an oracle for $c \in \CC$, such that $M_2^{\phi_1, \phi_2}$
halts if and only if $z$ converges to an attracting periodic orbit (or to $\infty$)
under $f_c: z \mapsto z^2+c$. 
\end{lem}

\begin{proof}
$M_2$ is systematically looking for an attracting cycle of $f_c$.
It also iterates $f_c$ on $z$ with increasing precision and for 
increasingly many steps until we are sure that either one of 
the two things holds:
\begin{enumerate}
\item 
the orbit of $z$ converges to $\infty$; or
\item 
we find an attracting orbit of $f_c$ and the orbit of $z$ converges to it. 
\end{enumerate}

If  the search is done systematically, the machine will eventually halt
if one of the possibilities above holds. It obviously won't halt if neither 
holds.
\end{proof}

\begin{proof} {\bf (of Theorem \ref{thm:main})} {\bf The algorithm is:} 
Run the machines $M_1^{\phi_1, \phi_2}(n)$ from Corollary \ref{cor:alg1} and $M_2^{\phi_1, \phi_2}$
from Lemma \ref{alg2} in parallel. Output $1$ if $M_1$ terminates first and $0$ if $M_2$ terminates
first. 

First we observe that $M_1(n)$ only halts on points that are $2^{-n}$-close to $\JJ$, in which 
case $1$ is a valid answer according to Definition \ref{weakdef}. Similarly, $M_2$ only halts
on points that are outside $\JJ$, in which case $0$ is a valid answer. Hence if the 
algorithm terminates, it outputs a valid answer. It remains to see that it does always terminate. 
Consider two cases.

{\bf Case 1: $(z,c) \in \overline{\JJ}$.} In this case $d((z,c),\JJ)=0<2^{-(n+4)}$, and 
the first machine is guaranteed to halt. 

{\bf Case 2: $(z,c) \notin \overline{\JJ}$.} By Lemma \ref{lem:onlyhyp}, $z$ converges to 
an attracting periodic orbit of $f_c$ in this case, and hence the second machine is 
guaranteed to halt. 
\end{proof}

\section{Proof of Lemma \ref{lem:onlyhyp}}
\label{sec:perturb}

Suppose $z\notin J_c$ and the orbit of $z$ does not belong to an attracting
basin. By the Fatou-Sullivan classification (see e.g. \cite{Mil}), 
there exists $k\in\NN$ such that $w\equiv f_c^k(z)$ belongs to a Siegel
disk or to the  immediate basin of a parabolic orbit.
Our aim is to show that for an arbitrary small $\delta>0$, there exists
a pair $(\tl z,\tl c)\in\CC\times\CC$ with $|z-\tl z|<\delta$,
$|c-\tl c|<\delta$, and for which $\tl z\in J_{\tl c}$.
We will treat the Siegel case first.

\subsection{The case when $w$ lies in a Siegel disk}
Let us denote $\Delta$ the Siegel disk containing $w$, and let $m\in\NN$ be its
period, that is, the mapping $$f_c^m:\Delta\to\Delta$$ is 
conjugated  by a conformal change of coordinates $\phi:\Delta\to \DD$
to an irrational rotation of $\DD$. 

The following statement is elementary (cf. Prop. 7.1 in \cite{Dou}):

\begin{prop}
\label{perturbsiegel}
Denote $\zeta=\phi^{-1}(0)\in\Delta$ the center of the Siegel disk.
For each $s>0$ there exists $\tl c\in B(c,s)$ such that
$f_{\tl c}$ has a parabolic periodic point $\tl \zeta$ of period $m$ 
in $B(\zeta,s)$. In particular,
$J_{\tl c}$ is connected, and 
 $B(\zeta,s)\cap J_{\tl c}\neq\emptyset.$
\end{prop}

Consider now the $f_c^m$-invariant analytic circle 
$$S_r=\phi^{-1}(\{z=re^{2\pi i\theta},\;\theta\in[0,2\pi)\})$$
which contains $w$. Let $\eps>0$ be such that 
$$B(w,\eps)\subset f_c^k(B(z,\delta))\cap \Delta.$$
Set $B\equiv B(w,\eps/2)$ and 
let $n\in\NN$ be such that the union
$$\bigcup_{0\leq i\leq n}f^{mi}_c(B)\supset S_r.$$
By Proposition \ref{perturbsiegel} for all $\delta>0$ small enough, there exist
$\tl c\in B(c,\delta)$ for which $J_{\tl c}$ is connected and 
there is a point of $J_{\tl c}$ inside
the domain bounded by $S_r$. Since repelling periodic orbits of $f_c$
are dense in $\partial \Delta$, again for $\delta$ small enough,
there are points of $J_{\tl c}$ on the outside of $S_r$ as well, and
so there exists a point $\xi\in J_{\tl c}\cap S_r$.
By construction, there exists $j\in\NN$ such that 
$f_c^j(B(z,\delta))\ni \xi$. By invariance of Julia set, if
$\tl c$ is close enough to $c$ we have 
$B(z,\delta)\cap J_{\tl c}\neq\emptyset$, and the proof is complete.

\subsection{The case when $w$ lies in a parabolic basin}
We need to recall the Douady-Lavaurs theory of parabolic implosion
(\cite{Dou,Lav}). Denote $\zeta$ the parabolic periodic point of  
$f_c$ whose immediate basin contains $w$, and let $m\in\NN$ be its
period. 

Recall that an attracting petal $P_A$ is a topological disk whose
boundary contains $\zeta$, such that $f_c^{mk}(P_A)\subset P_A$
for some $k\in\NN$, and such that the quotient Riemann surface
$$C_A=P_A/f_c^{mk}\simeq \CC/\ZZ.$$
The quotient $C_A$, called an attracting Fatou cylinder, parametrizes the
orbits converging under the dynamics of the
iterate $f_c^m$ to $\zeta$ inside the periodic cycle of petals
$P_A$, $f^m_c(P_A),\ldots,f^{m(k-1)}_c(P_A)$. 
Recall (see \cite{Mil}) that a quadratic polynomial
$f_c$ has only one cycle of petals at $\zeta$. A repelling petal $P_R$
is an attracting petal for the local inverse $f_c^{-m}$ fixing $\zeta$;
the union
$$\bigcup_{0\leq i\leq k-1}f_c^{mi}(P_A\cup P_R)$$
forms a neighborhood of $\zeta$. 
The repelling Fatou cylinder $C_R$ is defined in a similar fashion.

Let $\tau$ be any conformal isomorphism $C_A\to C_R$. After uniformization,
$$C_A\underset{\approx}{\mapsto}\CC/\ZZ,\;C_R\underset{\approx}{\mapsto}\CC/\ZZ$$
$\tau(z)\equiv z+q\mod\ZZ$ for some $q\in\CC$.
Let $g:P_A\to P_R$ be any lift of $\tau$; it necessarily commutes with $f_c^{mk}$. Consider the semigroup $G$ generated by the dynamics of the pair
$(f_c,g)$. The orbit $Gz$ of a point $z\in\CC$ is independent of the choice
of the lift $g$ and only depends on $\tau$.

Set $$J_{(c,\tau)}=\{z\in\CC\text{ such that }Gz\cap J_c\neq\emptyset\}.$$
It can be shown that this set is the boundary of 
$$K_{(c,\tau)}=\{z\in\CC\text{ such that }Gz\text{ is bounded}\}.$$

The Douady-Lavaurs theory postulates:

\begin{thm}
\label{explosion}
For every $\tau$ as above and every $s>0$ 
there exists $\tl c\in B(c,s)$ such that
$B(J_{\tl c},s)\supset J_{(c,\tau)}$.
\end{thm}

Since $\zeta\in J_c$, and $J_c$ is connected, there exists a point 
$u\in J_c\cap P_R$. Let $\hat w\in C_A$ be the orbit of $w$, and 
let $\hat u\in C_R$ be the orbit of $u$. Choose $\tau:C_A\to C_R$
so that $\tau(\hat w)=\hat u$. Then $J_{(c,\tau)}\ni z$, and the 
claim follows by Theorem \ref{explosion}.

 \end{document}